# CONVERGENCE RATES FOR BAYESIAN DENSITY ESTIMATION OF INFINITE-DIMENSIONAL EXPONENTIAL FAMILIES

By Catia Scricciolo

*University "L. Bocconi", Milan*

We study the rate of convergence of posterior distributions in density estimation problems for log-densities in periodic Sobolev classes characterized by a smoothness parameter $p$. The posterior expected density provides a nonparametric estimation procedure attaining the optimal minimax rate of convergence under Hellinger loss if the posterior distribution achieves the optimal rate over certain uniformity classes. A prior on the density class of interest is induced by a prior on the coefficients of the trigonometric series expansion of the log-density. We show that when $p$ is known, the posterior distribution of a Gaussian prior achieves the optimal rate provided the prior variances die off sufficiently rapidly. For a mixture of normal distributions, the mixing weights on the dimension of the exponential family are assumed to be bounded below by an exponentially decreasing sequence. To avoid the use of infinite bases, we develop priors that cut off the series at a sample-size-dependent truncation point. When the degree of smoothness is unknown, a finite mixture of normal priors indexed by the smoothness parameter, which is also assigned a prior, produces the best rate. A rate-adaptive estimator is derived.

**1. Introduction.** Bayesian nonparametrics is a very rapidly developing area of statistics. Several papers—including [1, 2, 4, 11, 14, 15, 16, 17, 18, 20, 21, 24, 25, 26]—have been devoted to the investigation of asymptotic properties of posterior distributions on infinite-dimensional parameter spaces.

The problem of estimating a density function $f_0$ w.r.t. the Lebesgue measure $\lambda$ on the unit interval, given a sample of i.i.d. observations $X_1, \ldots, X_n$ from $f_0$, is considered from the Bayesian perspective. Suppose that the sampling probability measure $P_0$ lies in $\mathscr{F}$, a class of absolutely continuous probability measures w.r.t. $\lambda$, equipped with the Hellinger metric $d_H$, the $L_2$-distance between square-rooted densities. Suppose, further, that the generic

density is of the form

$$\frac{\exp\{\theta(x)\}}{\int_0^1 \exp\{\theta(t)\}\,dt}, \qquad x \in [0,1],$$

with $\theta$ a square-integrable function in a periodic Sobolev class. We recall that for any given integer $p \geq 1$ and real $L > 0$, the Sobolev functional class $W(p, L)$ comprises all square-integrable functions with absolutely continuous derivative of order $p-1$ and $p$th derivative bounded in $L_2$-norm,

$$W(p, L) = \{\theta \in L_2[0,1] : \theta^{(p-1)} \text{ is abs. cont., } \|\theta^{(p)}\|_2^2 < L^2\}.$$

The periodic Sobolev class $W^{per}(p, L)$ is the following subclass of all periodic functions with period 1 satisfying the boundary conditions indicated:

$$W^{per}(p, L) = \{\theta \in W(p, L) : \theta^{(r)}(0) = \theta^{(r)}(1),\ r = 0, \ldots, p-1\}.$$

The problem is made discrete by representing a periodic Sobolev class as a Sobolev ellipsoid of $\ell_2$ via the trigonometric series expansion. Let $\{\phi_j(\cdot), j = 0, 1, \ldots\}$ be the orthonormal trigonometric system of $L_2[0,1]$. For $x \in [0,1]$, $\phi_0(x) \equiv 1$ and for $k \geq 1$, $\phi_{2k-1}(x) = \sqrt{2}\sin(2\pi kx)$, $\phi_{2k}(x) = \sqrt{2}\cos(2\pi kx)$. For $\theta \in L_2[0,1]$, let $\theta_j = \int_0^1 \theta(x)\phi_j(x)\,dx$, $j \geq 0$, be the sequence of its Fourier coefficients. To ease the notation, let $\boldsymbol{\theta} = (\theta_0, \theta_1, \ldots)$, $\boldsymbol{\phi}(\cdot) = (\phi_0(\cdot), \phi_1(\cdot), \ldots)$ and $\boldsymbol{\theta} \cdot \boldsymbol{\phi}(\cdot) = \sum_{j=0}^\infty \theta_j \phi_j(\cdot) = \theta_0 + \sum_{j=1}^\infty \theta_j \phi_j(\cdot)$. Each $\theta$ having the series expansion

$$\theta(x) = \boldsymbol{\theta} \cdot \boldsymbol{\phi}(x) = \theta_0 + \sqrt{2}\sum_{k=1}^\infty [\theta_{2k-1}\sin(2\pi kx) + \theta_{2k}\cos(2\pi kx)], \qquad x \in [0,1],$$

lies in $W^{per}(p, L)$ if and only if $\boldsymbol{\theta}$ belongs to the Sobolev ellipsoid of $\ell_2$,

$$E_p(Q) = \left\{\boldsymbol{\theta} \in \ell_2 : \sum_{j=0}^\infty v_j^{2p}\theta_j^2 < Q\right\}, \qquad Q = \frac{L^2}{\pi^{2p}},$$

with

$$v_0 = 0, \qquad v_j = \begin{cases} j+1, & \text{for } j \text{ odd,} \\ j, & \text{for } j \text{ even,} \end{cases} \qquad j = 1, 2, \ldots.$$

For $Q = \infty$, the Sobolev space $\{\boldsymbol{\theta} \in \ell_2 : \sum_{j=0}^\infty v_j^{2p}\theta_j^2 < \infty\}$ will be denoted by $E_p$. Setting $\psi(\boldsymbol{\theta}) = \log(\int_0^1 \exp\{\boldsymbol{\theta}\cdot\boldsymbol{\phi}(t)\}\,dt)$, the generic density can be rewritten as

$$f_{\boldsymbol{\theta}}(x) = \exp\{\boldsymbol{\theta}\cdot\boldsymbol{\phi}(x) - \psi(\boldsymbol{\theta})\} = \frac{\exp\{\sum_{j=1}^\infty \theta_j\phi_j(x)\}}{\int_0^1 \exp\{\sum_{j=1}^\infty \theta_j\phi_j(t)\}\,dt}, \qquad x \in [0,1].$$

The form of $f_{\boldsymbol{\theta}}$ explains why $\mathscr{F}$, which will also denote the density class $\{f_{\boldsymbol{\theta}},\ \boldsymbol{\theta} \in E_p(Q)\}$, is called an infinite-dimensional exponential family. Since $f_{\boldsymbol{\theta}}$ does not depend on $\theta_0$, for any pair $\boldsymbol{\theta}, \boldsymbol{\theta}' \in E_p(Q)$, the corresponding probability measures $P_{\boldsymbol{\theta}}, P_{\boldsymbol{\theta}'}$ are such that $P_{\boldsymbol{\theta}} \neq P_{\boldsymbol{\theta}'}$ if and only if $\theta_j \neq \theta'_j$ for some $j \geq 1$. Sequences differing only in the first coordinate identify the same probability measure and, thus, form an equivalence class. For instance, for $f_0 \in \mathscr{F}$, any $\boldsymbol{\theta}_0 \in E_p(Q)$ such that $f_{\boldsymbol{\theta}_0} = f_0$ can be taken as a representative of the class. It is now useful to highlight the fact that the $f_{\boldsymbol{\theta}}$'s are uniformly bounded and bounded away from zero. Let $\|\phi_j\|_\infty = \sup_{0 \leq x \leq 1} |\phi_j(x)|$ be the supremum norm of $\phi_j$, $j \geq 0$. Note that $\|\phi_j\|_\infty = \sqrt{2}$ for all $j \geq 1$. Setting

$$A = \sum_{j=1}^\infty v_j^{-2p}, \qquad B = \sqrt{2QA},$$

for each $\boldsymbol{\theta} \in E_p(Q)$, we have

$$\|\boldsymbol{\theta} \cdot \boldsymbol{\phi} - \theta_0\|_\infty \leq \sqrt{2} \sum_{j=1}^\infty |\theta_j| \leq \sqrt{2} \sqrt{\sum_{j=1}^\infty v_j^{2p} \theta_j^2} \sqrt{\sum_{j=1}^\infty v_j^{-2p}} < B < \infty.$$

Consequently, $\sup_{\boldsymbol{\theta} \in E_p(Q)} \|f_{\boldsymbol{\theta}}\|_\infty < e^{2B}$. Thus, the Hellinger distance between any pair $P_{\boldsymbol{\theta}'}, P_{\boldsymbol{\theta}} \in \mathscr{F}$, $d_{\mathrm{H}}(P_{\boldsymbol{\theta}'}, P_{\boldsymbol{\theta}}) = \{\int_0^1 (\sqrt{f_{\boldsymbol{\theta}'}} - \sqrt{f_{\boldsymbol{\theta}}})^2 \, d\lambda\}^{1/2}$, the Kullback–Leibler divergence $K(P_{\boldsymbol{\theta}'} \| P_{\boldsymbol{\theta}}) = K(f_{\boldsymbol{\theta}'} \| f_{\boldsymbol{\theta}}) = \int_0^1 f_{\boldsymbol{\theta}'} \log(f_{\boldsymbol{\theta}'}/f_{\boldsymbol{\theta}}) \, d\lambda$ and the $L_2$-distance $\|f_{\boldsymbol{\theta}'} - f_{\boldsymbol{\theta}}\|_2$ are equivalent and can be interchangeably used as loss functions.

The problem of estimating densities from exponential families has been studied by Crain [7, 8, 9, 10] from the frequentist perspective, where log-densities are generated by Legendre polynomials on $[-1, 1]$. Verdinelli and Wasserman [21] have used the same model for Bayesian goodness-of-fit testing. Our goal is to develop Bayesian density estimators attaining the optimal rate of convergence in the minimax sense under Hellinger loss, which is well known to be $n^{-p/(2p+1)}$ (see, for example, [25], Corollary 1, page 1574),

$$\inf_{\hat{f} \in \mathscr{S}_n} \sup_{f_{\boldsymbol{\theta}} \in \mathscr{F}} \mathbb{E}^n_{\boldsymbol{\theta}}[K(f_{\boldsymbol{\theta}} \| \hat{f})] \asymp \inf_{\hat{f} \in \mathscr{S}_n} \sup_{f_{\boldsymbol{\theta}} \in \mathscr{F}} \mathbb{E}^n_{\boldsymbol{\theta}}[d^2_{\mathrm{H}}(f_{\boldsymbol{\theta}}, \hat{f})]$$

$$\asymp \inf_{\hat{f} \in \mathscr{S}_n} \sup_{f_{\boldsymbol{\theta}} \in \mathscr{F}} \mathbb{E}^n_{\boldsymbol{\theta}}[\|f_{\boldsymbol{\theta}} - \hat{f}\|_2^2] \asymp n^{-2p/(2p+1)},$$

where $\mathscr{S}_n$ is the set of all estimators $\hat{f}$ for densities $f_{\boldsymbol{\theta}}$ in $\mathscr{F}$ based on $n$ observations and the expectation is taken over the $n$-fold product measure of $P_{\boldsymbol{\theta}}$. By writing $a_n \asymp b_n$, we mean that both $a_n \lesssim b_n$ and $b_n \lesssim a_n$, where $a_n \lesssim b_n$ if $a_n = O(b_n)$, namely, if there exists a constant $c$ such that $a_n \leq cb_n$ for all large $n$. Hereafter, all symbols $O$ and $o$ will refer to asymptotics as $n \to \infty$. The posterior expected density, which will be referred to as the *Bayes' estimator* and denoted by $\hat{f}_n$ in what follows, is a natural and common procedure for density estimation. From the general theory concerning posterior rates of convergence, it is known that if the posterior distribution on $\mathscr{F}$ converges at the exponential rate $e^{-Cn\varepsilon_n^2}$, where $\varepsilon_n$ is a positive sequence such that $\varepsilon_n \to 0$ and $n\varepsilon_n^2 \to \infty$ as $n \to \infty$, then the Bayes' estimator converges to the true density $f_0$ in the Hellinger distance at least as fast as $\varepsilon_n$ (see, e.g., [14], pages 506–507). Therefore, it suffices to put priors on $\mathscr{F}$ such that the corresponding posterior distributions converge exponentially fast at the optimal rate $n^{-p/(2p+1)}$. Recall that for $P_0 \in \mathscr{F}$, if $\Pi_n$ is a prior on $\mathscr{F}$ possibly depending on the sample size, the posterior converges at rate $\varepsilon_n$ (relative to $d_{\mathrm{H}}$) if for every positive sequence $M_n \to \infty$ such that $M_n \varepsilon_n \to 0$, $\Pi_n(H_{\varepsilon_n}^c(P_0)|X_1,\ldots,X_n) \to 0$ as $n \to \infty$, in probability or almost surely when sampling from $P_0$, where $H_{\varepsilon_n}^c(P_0) = \{P_{\boldsymbol{\theta}} \in \mathscr{F}: d_{\mathrm{H}}(P_0, P_{\boldsymbol{\theta}}) > M_n \varepsilon_n\}$. Since any prior on $E_p(Q)$ induces a prior on $\mathscr{F}$ via the map $\boldsymbol{\theta} \mapsto f_{\boldsymbol{\theta}}$, we can conveniently work with priors for the Fourier coefficients. Hereafter, we state a sufficient condition for posterior convergence at the optimal rate. The proof, deferred to the Appendix, relies on the fact that, in the present setting, Hellinger neighborhoods of $P_0$ translate into $\ell_2$-neighborhoods of $\boldsymbol{\theta}_0$.

THEOREM 1. *Let $\pi_n$ be a sequence of priors on $E_p(Q)$ and $\Pi_n$ the sequence of priors induced on $\mathscr{F}$. Suppose $\boldsymbol{\theta}_0 \in E_p(Q)$. Let $B_1^2 = e^{-8B}$ and $\varepsilon_n = n^{-p/(2p+1)}$. If for constants $c_1$, $c_2 > 0$,*

$$\pi_n\left(\left\{\boldsymbol{\theta}: \sum_{j=1}^{\infty}(\theta_j - \theta_{0,j})^2 \leq B_1^2 \varepsilon_n^2\right\}\right) \geq c_2 e^{-c_1 n \varepsilon_n^2}, \tag{1}$$

*then for a sufficiently large constant $M > 0$,*

$$\Pi_n(\{P_{\boldsymbol{\theta}}: d_{\mathrm{H}}(P_0, P_{\boldsymbol{\theta}}) > M\varepsilon_n\}|X_1,\ldots,X_n) \to 0 \quad \text{as } n \to \infty,$$

$P_0^{\infty}$-*almost surely, where $P_0^{\infty}$ denotes the infinite product measure of $P_0$.*

We develop several priors yielding Bayes' estimators that attain the optimal minimax rate. Preliminary ascertainment of consistency is based on results by Barron, Schervish and Wasserman [2], Walker and Hjort [23] and Walker [22], who have addressed the issue of consistency of posterior distributions for infinite-dimensional exponential families generated by orthonormal systems of bounded basis functions where the $\theta_j$'s are independent, zero-mean normals with variances chosen to ensure that $f_{\boldsymbol{\theta}}$ is a density with prior probability one. Then $K(P_0\|\lambda) < \infty$ is a sufficient condition for strong consistency.

We begin by considering the case where $p$ is known. In Section 2, we show that a sample-size-dependent prior constructed from an infinite product of normals achieves the optimal rate provided the variances decay sufficiently fast. The corresponding Bayes' estimator attains the minimax rate over Sobolev ellipsoids. As shown in Section 3, it is also attained by the

ESTIMATION OF EXPONENTIAL FAMILIES 5posterior expected density arising from a mixture of normals with mixing weights on the family dimension $k$ that are bounded below by a sequence exponentially decaying in $k$. Both estimators involve infinitely many basis functions. Thus, the need arises to develop priors on finite sets of coefficients. This implies truncating the series at a maximum number of components that is allowed to grow with sample size. Approximate density estimators are derived in Section 4. In Section 5, we consider the case where the degree of smoothness of $f_0$ is unknown. A prior on the smoothness parameter is assigned that has finite support. Normal priors with dimension depending on the smoothness parameter are combined into an overall distribution whose posterior is seen to converge at the best rate. An adaptive estimator is constructed. Adaptive convergence rates for posterior distributions on infinite-dimensional exponential families generated by wavelets with coefficients in a Besov space have been studied by Huang [16]. The relationship between our work and this article is considered in Section 6, along with some other closing remarks.

**2. Priors constructed from infinite normals.** A prior for $\boldsymbol{\theta}$ results from assuming independent, zero-mean normal coordinates. If we take $\theta_j \sim N(0, \tau_j^2)$, $j \geq 0$, with $\sum_{j=0}^{\infty} \tau_j^2 < \infty$, then the $\tau_j^2$'s must be specified so that the infinite product measure gives positive probability to $E_p$. Hereafter, we shall use $\lfloor x \rfloor$ ($\lceil x \rceil$) to mean the greatest (least) integer less (greater) than or equal to $x$. For each $n \geq 1$, let $\varepsilon_n = n^{-p/(2p+1)}$ and define $N_n = \lceil (8Q/(B_1^2 \varepsilon_n^2))^{1/(2p)} \rceil$, with $B_1^2 = e^{-8B}$ as before. We omit the subscript $n$ in $N_n$. Let $\tau_0^2 = 0$, which corresponds to a point mass at zero for the prior of $\theta_0$. Also, let $\tau_j^2 = \sigma^2 v_j^{-2q}$, with $q = p + 1/2$ for $j = 1, \ldots, N$, and $q = 2p + \alpha$, with $\alpha > 1/2$, for $j \geq N+1$. With this choice,

$$\sum_{j=0}^{\infty} v_j^{2p} \tau_j^2 = \sigma^2 \sum_{j=1}^{N} v_j^{2p} v_j^{-(2p+1)} + \sigma^2 \sum_{j=N+1}^{\infty} v_j^{2p} v_j^{-(4p+2\alpha)} < \infty,$$

hence, $\sum_{j=0}^{\infty} v_j^{2p} \theta_j^2$ converges almost surely; see (5.13) in [26], page 541. Let $\mu_n$ denote the sample-size-dependent prior

$$\mu_n(\boldsymbol{\theta}) = \boldsymbol{\delta}_0 \times \prod_{j=1}^{N} \frac{1}{\sigma v_j^{-(p+1/2)}} \phi\left(\frac{\theta_j}{\sigma v_j^{-(p+1/2)}}\right)$$
$$\times \prod_{j=N+1}^{\infty} \frac{1}{\sigma v_j^{-(2p+\alpha)}} \phi\left(\frac{\theta_j}{\sigma v_j^{-(2p+\alpha)}}\right), \qquad \boldsymbol{\theta} \in \mathbb{R}^{\infty},$$

where $\boldsymbol{\delta}_0$ denotes a point mass at zero, $\phi$ stands for the standard normal density and $\mathbb{R}^{\infty}$ is the space of sequences of real numbers. Let $\pi_n$ be the restriction of $\mu_n$ to $E_p(Q)$,

$$\pi_n(\boldsymbol{\theta}) = \frac{I_{E_p(Q)}(\boldsymbol{\theta})\mu_n(\boldsymbol{\theta})}{\mu_n(E_p(Q))}, \qquad \boldsymbol{\theta} \in \mathbb{R}^\infty. \tag{2}$$

We prove that the posterior of $\Pi_n$, the prior induced on $\mathscr{F}$ by $\pi_n$, converges at optimal rate. Henceforth, we may set $\sigma^2 = 1$ without loss of generality because the results of the following theorem and Corollary 1 are not affected by the value of $\sigma^2$ up to constants.

THEOREM 2. *If $\boldsymbol{\theta}_0 \in E_p(Q)$, then for a sufficiently large constant $M > 0$,*

$$\Pi_n(\{P_{\boldsymbol{\theta}} : d_{\mathrm{H}}(P_0, P_{\boldsymbol{\theta}}) > Mn^{-p/(2p+1)}\}|X_1, \ldots, X_n) \to 0 \qquad \text{as } n \to \infty,$$

*$P_0^\infty$-almost surely.*

PROOF. In virtue of Theorem 1, we only need to show that condition (1) is satisfied. Clearly,

$$J_n \triangleq \pi_n\left(\left\{\boldsymbol{\theta} : \sum_{j=1}^\infty (\theta_j - \theta_{0,j})^2 \leq B_1^2 \varepsilon_n^2\right\}\right)$$

$$\geq \mu_n\left(\left\{\boldsymbol{\theta} \in E_p(Q) : \sum_{j=1}^\infty (\theta_j - \theta_{0,j})^2 \leq B_1^2 \varepsilon_n^2\right\}\right).$$

We show that for all large $n$,

$$J_n \geq \mu_n(E_n), \tag{3}$$

where

$$E_n = \left\{\boldsymbol{\theta} : \sum_{j=1}^N (\theta_j - \theta_{0,j})^2 \leq \frac{B_1^2 \varepsilon_n^2}{C_0}, \sum_{j=N+1}^\infty v_j^{2p} \theta_j^2 \leq \frac{B_1^2 \varepsilon_n^2}{8}\right\}$$

with $C_0$ a positive constant depending on $\boldsymbol{\theta}_0$ to be suitably chosen as will be prescribed. To prove (3), it suffices to show that for each $\boldsymbol{\theta} \in E_n$,

(i) $\boldsymbol{\theta} \in E_p(Q)$;
(ii) $\sum_{j=1}^\infty (\theta_j - \theta_{0,j})^2 < B_1^2 \varepsilon_n^2$.

We start with (i). Let $0 < \delta_0 \leq Q$ be such that $\sum_{j=0}^\infty v_j^{2p} \theta_{0,j}^2 = Q - \delta_0$. By Schwarz's inequality,

$$\sum_{j=0}^\infty v_j^{2p} \theta_j^2 \leq \sum_{j=1}^N v_j^{2p}(\theta_j - \theta_{0,j})^2 + \sum_{j=1}^N v_j^{2p} \theta_{0,j}^2$$

$$+ 2\sqrt{\sum_{j=1}^N v_j^{2p}(\theta_j - \theta_{0,j})^2} \sqrt{\sum_{j=1}^N v_j^{2p} \theta_{0,j}^2} + \sum_{j=N+1}^\infty v_j^{2p} \theta_j^2$$

$$\leq (N+1)^{2p}\frac{B_1^2\varepsilon_n^2}{C_0} + (Q-\delta_0)$$
$$+ 2\sqrt{(N+1)^{2p}\frac{B_1^2\varepsilon_n^2}{C_0}(Q-\delta_0)} + \frac{B_1^2\varepsilon_n^2}{8}.$$

Note that if $x > 0$, then for $0 < K \leq x$,

$$(\lceil x \rceil + 1)^{2p} \leq \lceil x \rceil^{2p}\left(1+\frac{1}{K}\right)^{2p} \leq (x+1)^{2p}\left(1+\frac{1}{K}\right)^{2p} \leq x^{2p}\left(1+\frac{1}{K}\right)^{4p}.$$

Fix $K \geq 1$ and let $n_1$ be the smallest $n$ such that $1 \leq K \leq (8Q/(B_1^2\varepsilon_n^2))^{1/(2p)}$. For $n \geq n_1$,

$$(N+1)^{2p} \leq \frac{8Q}{B_1^2\varepsilon_n^2}\left(1+\frac{1}{K}\right)^{4p} \leq 16^p \frac{8Q}{B_1^2\varepsilon_n^2}.$$

Fix $0 < \eta_0 \leq (\sqrt{Q} - \sqrt{Q-\delta_0})$ and define $C_0 = 16^{p+1}Q/\eta_0^2$. Let $n_2$ be the smallest $n$ such that $B_1^2\varepsilon_n^2/8 < \eta_0^2/2$. Obviously, $n_2$ depends on $\eta_0$. For each $n \geq \bar{n} = \max\{n_1, n_2\}$,

$$\sum_{j=0}^{\infty} v_j^{2p}\theta_j^2 < \frac{\eta_0^2}{2} + (Q-\delta_0) + 2\eta_0\sqrt{\frac{Q-\delta_0}{2}} + \frac{\eta_0^2}{2} \leq (\eta_0 + \sqrt{Q-\delta_0})^2 \leq Q,$$

which proves (i). We now turn to (ii). Using the inequality $(a+b)^2 \leq 2(a^2 + b^2)$, since $C_0 > 2$, we have

$$\sum_{j=1}^{\infty}(\theta_j - \theta_{0,j})^2 < \frac{B_1^2\varepsilon_n^2}{C_0} + 2\sum_{j=N+1}^{\infty} v_j^{2p}\theta_j^2 + 2N^{-2p}\sum_{j=N+1}^{\infty} v_j^{2p}\theta_{0,j}^2$$
$$< \frac{3B_1^2\varepsilon_n^2}{4} + 2N^{-2p}Q \leq B_1^2\varepsilon_n^2.$$

Hence, both (i) and (ii) are satisfied for all $n \geq \bar{n}$. We now find a lower bound on $\mu_n(E_n)$. By independence of the $\theta_j$'s,

$$J_n \geq \Pr\left(\left\{(\theta_1,\ldots,\theta_N): \sum_{j=1}^{N}(\theta_j - \theta_{0,j})^2 \leq \frac{B_1^2\varepsilon_n^2}{C_0}\right\}\right)$$
$$\times \Pr\left(\left\{(\theta_{N+1}, \theta_{N+2}, \ldots): \sum_{j=N+1}^{\infty} v_j^{2p}\theta_j^2 \leq \frac{B_1^2\varepsilon_n^2}{8}\right\}\right)$$
$$\triangleq J_{1,n} \times J_{2,n}.$$

Reasoning as in Lemma 4 of Shen and Wasserman [19], page 711, we obtain that

(4) $\quad J_{1,n} > e^{-(2Q+p+1/2)N} 2^{-(p+1)N} \Pr\left(\sum_{j=1}^{N} V_j^2 \leq 2\frac{B_1^2\varepsilon_n^2}{C_0}(2N)^{2p+1}\right),$

where $V_1, \ldots, V_N$ are independent, standard normal random variables. The probability on the right-hand side of (4) can be bounded below using Stirling's approximation. For ease of notation, let $\xi_n^2 = B_1^2 \varepsilon_n^2 / C_0$ and $d = p + 1/2$. Then

$$\Pr\left(\sum_{j=1}^{N} V_j^2 \leq 2(2N)^{2d} \xi_n^2\right) \gtrsim \frac{e^{-(2N)^{2d}\xi_n^2}(2N)^{dN}\xi_n^N}{(N/2)^{N/2-1}e^{-N/2}\sqrt{\pi N}}.$$

Noting that $(2N)^{2d}\xi_n^2 \leq (16^{p+1}Q/C_0)N = \eta_0^2 N$, we obtain that

$$J_{1,n} \gtrsim e^{-cN},$$

where $c = 2Q + p + \eta_0^2 - \frac{1}{2}\log(\eta_0^2/2^{4p+1}) > 0$. Let us consider $J_{2,n}$. By Markov's inequality,

$$J_{2,n} \geq 1 - \frac{8}{B_1^2 \varepsilon_n^2} \sum_{j=N+1}^{\infty} v_j^{2p} \mathbb{E}[\theta_j^2]$$

$$\geq 1 - \frac{8}{B_1^2 \varepsilon_n^2} \sum_{j=N+1}^{\infty} j^{-(2p+2\alpha)} > 1 - \frac{8Q/2}{B_1^2 \varepsilon_n^2 N^{2p}} \geq \frac{1}{2}$$

for all large $n$. Combining lower bounds on $J_{1,n}$ and $J_{2,n}$, we obtain that for $c_1 = 2c(8Q/B_1^2)^{1/(2p)}$ and all large $n$,

$$J_n \geq J_{1,n} \times J_{2,n} \gtrsim e^{-c_1 n \varepsilon_n^2},$$

which completes the proof. $\square$

COROLLARY 1. *If $\hat{f}_n$ is the Bayes' estimator arising from prior* (2), *then for any $0 < Q' < Q$,*

$$\sup_{\boldsymbol{\theta}_0 \in E_p(Q')} \mathbb{E}_{\boldsymbol{\theta}_0}^n[d_{\mathrm{H}}^2(f_{\boldsymbol{\theta}_0}, \hat{f}_n)] \asymp n^{-2p/(2p+1)}.$$

PROOF. Note that for each $\boldsymbol{\theta}_0 \in E_p(Q')$, choosing $0 < \eta \leq (\sqrt{Q} - \sqrt{Q'})$, Theorem 2 applies with constants that do not depend on the specific point $\boldsymbol{\theta}_0$. Thus, as a byproduct of Theorem A.1, for suitable constants $M$, $C$, $c > 0$ and sufficiently large $n$,

$$\sup_{\boldsymbol{\theta}_0 \in E_p(Q')} \mathbb{E}_{\boldsymbol{\theta}_0}^n[\Pi_n(\{P_{\boldsymbol{\theta}} : d_{\mathrm{H}}(P_{\boldsymbol{\theta}_0}, P_{\boldsymbol{\theta}}) > M\varepsilon_n\}|X_1, \ldots, X_n)] \leq Ce^{-cn\varepsilon_n^2}.$$

By Theorem 5 of Shen and Wasserman [19], page 694,

$$\sup_{\boldsymbol{\theta}_0 \in E_p(Q')} \mathbb{E}_{\boldsymbol{\theta}_0}^n[d_{\mathrm{H}}^2(f_{\boldsymbol{\theta}_0}, \hat{f}_n)] \leq M^2 \varepsilon_n^2 + 2Ce^{-cn\varepsilon_n^2} \lesssim \varepsilon_n^2,$$

which, combined with

$$\varepsilon_n^2 \lesssim \sup_{\boldsymbol{\theta}_0 \in E_p(Q')} \mathbb{E}_{\boldsymbol{\theta}_0}^n [d_{\mathrm{H}}^2(f_{\boldsymbol{\theta}_0}, \hat{f}_n)],$$

yields the assertion. □

REMARK 1. Corollary 1 shows that prior (2) yields a Bayes' density estimator attaining optimal minimax rate over any ellipsoid $E_p(Q')$, with $Q' < Q$. Theorem 2 and Corollary 1 are of interest because they establish that, for the problem under consideration, in contrast to the infinitely many normal means problem considered in [26], a sample-size-dependent direct Gaussian prior yields a Bayes' estimator attaining optimal minimax rate provided the prior variances die off sufficiently rapidly.

**3. Sieve priors.** In this section, we consider sieve priors restricted to $E_p(Q)$. Sieve priors have been used by Zhao [26] and Shen and Wasserman [19]. The basic idea is to put a prior on the dimension of the exponential family, hereafter denoted by $k$. Before describing the hierarchical structure of a sieve prior, we introduce some more notation. Henceforth, for any integer $N \geq 1$, let $\boldsymbol{\theta}_N = (\theta_0, \ldots, \theta_N, 0, 0, \ldots)$ denote a sequence such that all but possibly the first $N+1$ coordinates are equal to zero. Also, let $E_{p,N}(Q) = \{\boldsymbol{\theta}_N : \sum_{j=0}^N v_j^{2p} \theta_j^2 < Q\}$. Clearly, $E_{p,N}(Q) \subseteq E_p(Q)$.

(i) Conditionally on $k \geq 1$ and $\boldsymbol{\theta}$, for each $n \geq 1$, the random variables $X_1, \ldots, X_n$ are i.i.d., with density

$$f_{\boldsymbol{\theta}}(x) = \frac{\exp\{\sum_{j=1}^k \theta_j \phi_j(x)\}}{\int_0^1 \exp\{\sum_{j=1}^k \theta_j \phi_j(s)\} ds}, \qquad x \in [0,1];$$

(ii) conditionally on $k$, the sequence $\boldsymbol{\theta}$ has distribution $\mu_k$, which makes the coordinates independent and such that $\theta_0 \equiv 0$, $\theta_j \sim N(0, v_j^{-(2p+1)})$, $j = 1, \ldots, k$, and $\theta_j$ is degenerate at 0 for all $j > k$;
(iii) the exponential family dimension $k$ has distribution $\{\lambda(k), k = 1, 2, \ldots\}$ with $\lambda(k) \geq Ae^{-\gamma k}$, $k \geq 1$, for some $A, \gamma > 0$.

Let $\pi$ denote the restriction of the sieve prior $\mu = \sum_{k=1}^\infty \lambda(k) \mu_k$ to $E_p(Q)$,

(5) $$\pi(\boldsymbol{\theta}) = \frac{I_{E_p(Q)}(\boldsymbol{\theta}) \mu(\boldsymbol{\theta})}{\mu(E_p(Q))}, \qquad \boldsymbol{\theta} \in \mathbb{R}^\infty,$$

where $\mu(E_p(Q)) = \sum_{k=1}^\infty \lambda(k) \mu_k(E_{p,k}(Q))$. Next, we study the convergence rate for the posterior of the prior $\Pi$ induced by $\pi$ on $\mathscr{F}$.

THEOREM 3. *If $\boldsymbol{\theta}_0 \in E_p(Q)$, then for a sufficiently large constant $M > 0$,*

$$\Pi(\{P_{\boldsymbol{\theta}} : d_{\mathrm{H}}(P_0, P_{\boldsymbol{\theta}}) > Mn^{-p/(2p+1)}\}|X_1, \ldots, X_n) \to 0 \qquad as\ n \to \infty,$$

*$P_0^\infty$-almost surely.*

PROOF. We appeal to Theorem 1. Note that for $N = \lceil (2Q/(B_1^2 \varepsilon_n^2))^{1/(2p)} \rceil$,

(a) $$\sum_{j=N+1}^{\infty} \theta_{0,j}^2 < N^{-2p} \sum_{j=N+1}^{\infty} v_j^{2p} \theta_{0,j}^2 < N^{-2p} Q \leq B_1^2 \varepsilon_n^2/2.$$

Let $0 < \delta_0 \leq Q$ be such that $\sum_{j=0}^{\infty} v_j^{2p} \theta_{0,j}^2 = Q - \delta_0$. If $n$ is sufficiently large that $(2Q/(B_1^2 \varepsilon_n^2))^{1/(2p)} \geq 1$ and $\bar{B}_1$ is a positive constant such that $\bar{B}_1^2 < B_1^2 (1 - \sqrt{1 - \delta_0/Q})^2/2^{4p+1}$, then

(b) $$\left\{ \boldsymbol{\theta}_N \in E_{p,N}(Q) : \sum_{j=1}^{N} (\theta_j - \theta_{0,j})^2 \leq B_1^2 \varepsilon_n^2/2 \right\}$$
$$\supseteq \left\{ \boldsymbol{\theta}_N : \sum_{j=1}^{N} (\theta_j - \theta_{0,j})^2 \leq \bar{B}_1^2 \varepsilon_n^2 \right\}.$$

Using (a) and (b),

$$I_n \triangleq \pi\left( \left\{ \boldsymbol{\theta} : \sum_{j=1}^{\infty} (\theta_j - \theta_{0,j})^2 \leq B_1^2 \varepsilon_n^2 \right\} \right)$$
$$> \lambda(N) \mu_N \left( \left\{ \boldsymbol{\theta}_N \in E_{p,N}(Q) : \sum_{j=1}^{N} (\theta_j - \theta_{0,j})^2 + \sum_{j=N+1}^{\infty} \theta_{0,j}^2 \leq B_1^2 \varepsilon_n^2 \right\} \right)$$
$$\geq \lambda(N) \mu_N \left( \left\{ \boldsymbol{\theta}_N \in E_{p,N}(Q) : \sum_{j=1}^{N} (\theta_j - \theta_{0,j})^2 \leq B_1^2 \varepsilon_n^2/2 \right\} \right)$$
$$\geq \lambda(N) \mu_N \left( \left\{ \boldsymbol{\theta}_N : \sum_{j=1}^{N} (\theta_j - \theta_{0,j})^2 \leq \bar{B}_1^2 \varepsilon_n^2 \right\} \right) \triangleq \lambda(N) I_{1,n}.$$

Let $\zeta_n^2 = \bar{B}_1^2 \varepsilon_n^2$ and $d = p + 1/2$. Noting that $(2N)^{2d} \zeta_n^2 \leq (16^d Q \bar{B}_1^2 / B_1^2) N$ and $\sum_{j=1}^{N} v_j^{2d} \theta_{0,j}^2 < 2QN$, by Lemma 4 of Shen and Wasserman [19], page 711, and using Stirling's approximation, we obtain that

$$I_{1,n} > e^{-(2Q+d)N} 2^{-(d+1/2)N} \frac{1}{\Gamma(N/2)} \int_0^{(2N)^{2d} \zeta_n^2} z^{N/2-1} e^{-z} dz$$
$$\gtrsim e^{-(2Q+d)N} 2^{-(d+1/2)N} \frac{e^{-(2N)^{2d} \zeta_n^2} (2N)^{dN} \zeta_n^N}{(N/2)^{N/2-1} e^{-N/2} \sqrt{\pi N}} \gtrsim e^{-cN},$$

where $c = 2Q + p + 16^d Q \bar{B}_1^2 / B_1^2 - \frac{1}{2} \log(2Q \bar{B}_1^2 / B_1^2) > 0$. Therefore,

(6) $\quad I_n > \lambda(N) I_{1,n} \gtrsim e^{-(\gamma+c)N} \geq e^{-2(\gamma+c)(2Q/B_1^2)^{1/(2p)} \varepsilon_n^{-1/p}} = e^{-c_1 n \varepsilon_n^2},$

with $c_1 = 2(\gamma + c)(2Q/B_1^2)^{1/(2p)}$, and condition (1) is satisfied. $\square$

REMARK 2. An examination of the proof of Theorem 3 reveals that posterior convergence at the optimal rate depends on the assumed tail behavior of the mixing weights, which are bounded below by an exponentially decreasing sequence. This requirement is used in (6) to guarantee that $\varepsilon_n$-Hellinger-type neighborhoods of $P_0$ have prior mass at least of the order of $e^{-c_1 n \varepsilon_n^2}$.

COROLLARY 2. If $\hat{f}_n$ is the Bayes' estimator arising from prior (5), then for any $0 < Q' < Q$,
$$\sup_{\boldsymbol{\theta}_0 \in E_p(Q')} \mathbb{E}^n_{\boldsymbol{\theta}_0}[d^2_{\mathrm{H}}(f_{\boldsymbol{\theta}_0}, \hat{f}_n)] \asymp n^{-2p/(2p+1)}.$$

PROOF. It suffices to check that the convergence of the posterior is uniform over $E_p(Q')$. More formally, for each $\boldsymbol{\theta}_0 \in E_p(Q')$, Theorem 3 applies, with constants depending only on $Q$ and $Q'$, so that for suitable $M$, $C$, $c > 0$,
$$\sup_{\boldsymbol{\theta}_0 \in E_p(Q')} \mathbb{E}^n_{\boldsymbol{\theta}_0}[\Pi(\{P_{\boldsymbol{\theta}} : d_{\mathrm{H}}(P_{\boldsymbol{\theta}_0}, P_{\boldsymbol{\theta}}) > M\varepsilon_n\} | X_1, \ldots, X_n)] \leq C e^{-c n \varepsilon_n^2}.$$

Note that for $\delta' = Q - Q' > 0$, we have $\sum_{j=0}^{\infty} v_j^{2p} \theta_{0,j}^2 < Q' = Q - \delta'$ for all $\boldsymbol{\theta}_0 \in E_p(Q')$. Thus, Theorem 3 applies with $\bar{B}_1^2 < B_1^2(1 - \sqrt{1 - \delta'/Q})^2 / 2^{4p+1}$. The assertion then follows via reasoning similar to that used in the proof of Corollary 1. □

REMARK 3. Corollary 2 demonstrates that the Bayes' estimator attains the minimax rate of convergence under Hellinger loss over any ellipsoid $E_p(Q')$, with $Q' < Q$.

**4. Sample-size-dependent priors and density estimators.** Bayes' estimators arising from priors (2) and (5) involve infinitely many terms. To avoid the use of infinite bases, we define priors supported on exponential families whose dimension varies with sample size at a carefully chosen rate. Let $N_n$ be a sequence of positive integers, to be specified below. To ease the notation, we omit the subscript $n$ in $N_n$. For each $n \geq 1$, let $\mu_N$ be the prior that makes the coordinates independent and such that $\theta_0 \equiv 0$, $\theta_j \sim N(0, v_j^{-(2p+1)})$, $j = 1, \ldots, N$, and $\theta_j$ is degenerate at 0 for all $j > N$. Let

(7) $$\pi_N(\boldsymbol{\theta}_N) = \frac{I_{E_{p,N}(Q)}(\boldsymbol{\theta}_N) \mu_N(\boldsymbol{\theta}_N)}{\mu_N(E_{p,N}(Q))}, \qquad \boldsymbol{\theta}_N \in \mathbb{R}^{\infty},$$

be the restriction of $\mu_N$ to $E_{p,N}(Q)$ and let $\Pi_n$ denote the induced prior on $\mathscr{F}_n = \{f_{\boldsymbol{\theta}_N}, \boldsymbol{\theta}_N \in E_{p,N}(Q)\}$, where $f_{\boldsymbol{\theta}_N} = e^{\boldsymbol{\theta}_N \cdot \boldsymbol{\phi} - \psi(\boldsymbol{\theta}_N)}$.

THEOREM 4. *Let $N = \lceil (2Q/B_1^2)^{1/(2p)} n^{1/(2p+1)} \rceil$. If $\boldsymbol{\theta}_0 \in E_p(Q)$, then for a sufficiently large constant $M > 0$,*

$$\Pi_n(\{P_{\boldsymbol{\theta}} : d_H(P_0, P_{\boldsymbol{\theta}}) > Mn^{-p/(2p+1)}\} | X_1, \ldots, X_n) \to 0 \qquad as\ n \to \infty,$$

$P_0^\infty$*-almost surely.*

PROOF. The proof of Theorem 3 carries over to this case with simple modifications. □

REMARK 4. The assertion of Theorem 4 also holds true for the truncated sieve prior

$$\pi_n(\boldsymbol{\theta}) = \frac{I_{E_{p,N}(Q)}(\boldsymbol{\theta}) \mu_n(\boldsymbol{\theta})}{\mu_n(E_{p,N}(Q))}, \qquad \boldsymbol{\theta} \in \mathbb{R}^\infty, \tag{8}$$

where for each $n \geq 1$, $\mu_n = \sum_{k=1}^N \lambda_n(k) \mu_k$, with $\lambda_n(k) \geq A_1 e^{-\gamma k}$, $k = 1, \ldots, N$, and $\sum_{k=1}^N \lambda_n(k) = 1$. A uniform version of Theorem 4 can be formulated for priors (7) and (8) so that the corresponding Bayes' estimators attain minimax rate.

In the next proposition, approximations for the Bayes' estimators arising from priors (7) and (8) are provided.

PROPOSITION 1. *If for given (large) $n$, $Q$ is such that $\mu_N(E_{p,N}(Q)) \simeq 1$ ($\mu_n(E_{p,N}(Q)) \simeq 1$), then the Bayes' estimators arising from priors (7) and (8) can be approximated by*

$$C_{1,n} \exp\left\{ \frac{1}{2} \sum_{j=1}^N \frac{\phi_j^2(x) + 2n \bar{\phi}_j \phi_j(x)}{v_j^{2p+1} + n + 1} \right\}, \qquad x \in [0, 1], \tag{9}$$

*and*

$$C_{2,n} \sum_{k=1}^N \lambda_n(k) \rho_n(k) \exp\left\{ \frac{1}{2} \sum_{j=1}^k \frac{(\phi_j(x) + n \bar{\phi}_j)^2}{v_j^{2p+1} + n + 1} \right\}, \qquad x \in [0, 1], \tag{10}$$

*respectively, where $N$ is defined as in Theorem 4, $\bar{\phi}_j = n^{-1} \sum_{i=1}^n \phi_j(X_i)$, $j = 1, \ldots, N$, $\rho_n(k) = \prod_{j=1}^k (1 + (n+1) v_j^{-(2p+1)})^{-1/2}$, $k = 1, \ldots, N$, and $C_{1,n}$, $C_{2,n}$ stand for the normalizing constants.*

PROOF. First, note that for given $n$, if $Q$ is sufficiently large, then $\mu_N(E_{p,N}(Q)) \simeq 1$. To see this, observe that since $\theta_0$ is degenerate at zero, the probability $\mu_N(E_{p,N}(Q))$ is bounded below by the left tail of the chi-square distribution with $N$ degrees of freedom,

$$\mu_N(E_{p,N}(Q)) \geq \mu_N\left(\left\{\boldsymbol{\theta}_N : \sum_{j=1}^N v_j^{2p+1} \theta_j^2 < Q\right\}\right) = \Pr(\chi_N^2 < Q).$$

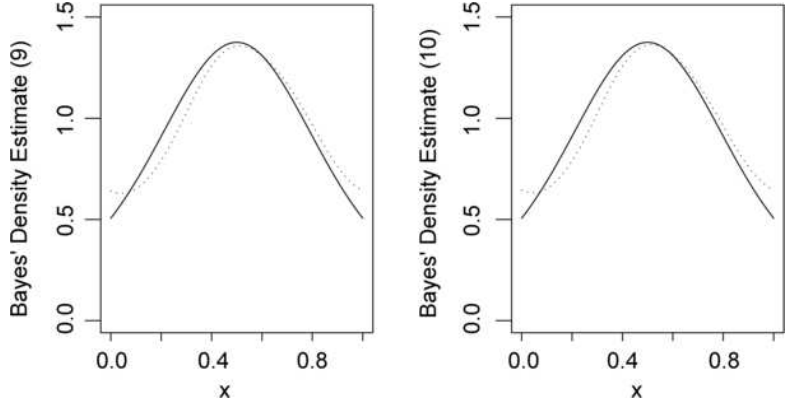

Fig. 1. *True density (solid line) and its approximate Bayes' estimate* (9) *(dotted line) on the left. True density (solid line) and its approximate Bayes' estimate* (10) *(dotted line) on the right.*

Similarly, $\mu_n(E_{p,N}(Q)) \geq \sum_{k=1}^{N} \lambda_n(k) \Pr(\chi_k^2 < Q) \geq \Pr(\chi_N^2 < Q)$ because the chi-square distribution is stochastically increasing in its degrees of freedom. We now derive (9). Setting $L_n(f_{\boldsymbol{\theta}_N}) = \prod_{i=1}^{n} f_{\boldsymbol{\theta}_N}(X_i)$, the posterior expected density can be written as

$$\hat{f}_n(x) = \frac{\mathbb{E}[f_{\boldsymbol{\theta}_N}(x) L_n(f_{\boldsymbol{\theta}_N})]}{\int_0^1 \mathbb{E}[f_{\boldsymbol{\theta}_N}(s) L_n(f_{\boldsymbol{\theta}_N})] \, ds}, \qquad x \in [0,1],$$

where $\mathbb{E}$ stands for expectation under prior (7). Since $\mu_N(E_{p,N}(Q)) \backsimeq 1$, $\mu_N$ can be thought of as a prior on $\mathbb{R}^N$ and

$$\mathbb{E}[f_{\boldsymbol{\theta}_N}(x) L_n(f_{\boldsymbol{\theta}_N})] \backsimeq \int_{\mathbb{R}^N} f_{\boldsymbol{\theta}_N}(x) L_n(f_{\boldsymbol{\theta}_N}) \mu_N(d\boldsymbol{\theta}_N), \qquad x \in [0,1].$$

Since $n$ is large, we can proceed as in Corollary 1 of Lenk [17], pages 534–535 (see also pages 541–542), and approximate $e^{(n+1)\psi(\boldsymbol{\theta}_N)}$ using the CLT. Straightforward computations then lead to (9). Approximation (10) may be proved similarly. □

REMARK 5. The number of terms $N = O(n^{1/(2p+1)})$ used in (9) and (10) is of the same order as the dimension, say $N^*$, of the exponential family employed to define the density estimator proposed by Barron and Sheu [3], when the log-density is in the periodic Sobolev space $W^{per}(p, \infty)$. Such an estimator, say $\hat{f}$, is defined to maximize the likelihood in the $N^*$-dimensional exponential family and is shown to converge to $f_0$ in the sense of relative entropy (Kullback–Leibler divergence) at rate $O_P(n^{-2p/(2p+1)})$, that is, $K(f_0 \| \hat{f}) = O_P(n^{-2p/(2p+1)})$.

The plots in Figure 1 show approximate Bayes' estimates (9) on the left-hand side and (10) on the right-hand side for the density function

$$\frac{\exp\{\sin(\pi x)\}}{\int_0^1 \exp\{\sin(\pi t)\}\, dt}, \qquad x \in [0,1],$$

based on $n = 500$ observations. We took $p = 2$, $N = O(n^{1/5})$ and $\lambda_n(k) \propto e^{-\gamma k}$, with $\gamma = 0.1$. Both estimates, which appear very similar, are close to the true density.

**5. Rate adaptation.** Thus far, we have assumed that the degree of smoothness, $p$, of $f_0$ is known. We now suppose that this is unknown and denote its value by $p_0$. In accordance with the Bayesian approach, we may consider $p$ as a hyperparameter and assign it a prior distribution. Let $\mathbb{P} = \{p_{\underline{m}}, \ldots, p_{-1}, p_0, p_1, \ldots, p_{\overline{m}}\}$ be a finite set of possible values for $p$, with $1 \leq p_{\underline{m}} < \cdots < p_{-1} < p_0 < p_1 < \cdots < p_{\overline{m}}$. Let $\mathbb{M} = \{\underline{m}, \ldots, -1, 0, 1, \ldots, \overline{m}\}$ be the corresponding index set. For any $m \in \mathbb{M}$, let $N_m = \lceil n^{1/(2p_m+1)} \rceil$, where the subscript $m$ is introduced to stress the dependence on $p_m$. We consider the following hierarchical prior. For each $n \geq 1$,

(i) conditionally on $p = p_m$ and $\boldsymbol{\theta}$, the random variables $X_1, \ldots, X_n$ are i.i.d. with density

$$f_{\boldsymbol{\theta}}(x) = \frac{\exp\{\sum_{j=1}^{N_m} \theta_j \phi_j(x)\}}{\int_0^1 \exp\{\sum_{j=1}^{N_m} \theta_j \phi_j(s)\}\, ds}, \qquad x \in [0,1];$$

(ii) conditionally on $p = p_m$, $\boldsymbol{\theta}$ has distribution $\mu_{N_m}$, which makes the coordinates independent and such that $\theta_0 \equiv 0$, $\theta_j \sim N(0, v_j^{-(2p_m+1)})$, $j = 1, \ldots, N_m$, and $\theta_j$ is degenerate at 0 for all $j > N_m$;
(iii) $p$ has distribution $w(m) = \Pr(p = p_m) > 0$ for all $m \in \mathbb{M}$.

The overall prior is $\pi_n = \sum_{m \in \mathbb{M}} w(m) \mu_{N_m}$. Let $\Pi_n$ be the induced prior on $\bigcup_{m \in \mathbb{M}} \{f_{\boldsymbol{\theta}_{N_m}}, \boldsymbol{\theta}_{N_m} \in \mathbb{R}^\infty\}$. Our goal is to show that this mixture prior achieves the rate of convergence $n^{-p_0/(2p_0+1)}$ whenever $\boldsymbol{\theta}_0 \in E_{p_0}$, with $p_0 \in \mathbb{P}$. We need to introduce further notation. For each $j \geq 1$, let $\mathbb{E}_0[\phi_j(X_1)]$ and $\mathbb{V}_0[\phi_j(X_1)]$ be the expected value and variance of $\phi_j(X_1)$ w.r.t. $P_0$, respectively. Note that $\mathbb{E}_0[\phi_j(X_1)] \leq \sqrt{2}$ and $\mathbb{V}_0[\phi_j(X_1)] \leq 2$ for all $j \geq 1$. The conditions

$$\sum_{j=1}^{\infty} v_j^{2p_0} (\mathbb{E}_0[\phi_j(X_1)])^2 < \infty, \tag{11}$$

$$\sum_{j=1}^{\infty} \mathbb{V}_0[\phi_j(X_1)] < \infty \tag{12}$$

are assumed to be in force in what follows. We are now in a position to state the main result of this section.

THEOREM 5. *Suppose $p_0 \in \mathbb{P}$. If $\boldsymbol{\theta}_0 \in E_{p_0}$ satisfies conditions (11) and (12), then for a sufficiently large constant $M > 0$,*

$$\Pi_n(\{P_{\boldsymbol{\theta}} : d_{\mathrm{H}}(P_0, P_{\boldsymbol{\theta}}) > M n^{-p_0/(2p_0+1)}\} | X_1, \ldots, X_n) \to 0$$

*in $P_0^n$-probability as $n \to \infty$.*

PROOF. The idea is to show that the posterior mass will ultimately be lying in a Sobolev ellipsoid. This will drastically reduce the effective parameter space, allowing us to apply the theory developed above. Let $\varepsilon_n = n^{-p_0/(2p_0+1)}$. Define $\overline{w}(m) = w(m)/\sum_{l \geq 0} w(l)$, for $m = 0, \ldots, \overline{m}$, and let $\overline{\pi}_n = \sum_{m \geq 0} \overline{w}(m) \mu_{N_m}$. For any $Q > 0$,

$$U_n \triangleq \Pi_n(\{P_{\boldsymbol{\theta}} : d_{\mathrm{H}}(P_0, P_{\boldsymbol{\theta}}) > M\varepsilon_n\} | X_1, \ldots, X_n)$$
$$= \pi_n(\{\boldsymbol{\theta} : d_{\mathrm{H}}(P_0, P_{\boldsymbol{\theta}}) > M\varepsilon_n\} | X_1, \ldots, X_n)$$
$$= \Pr(\{\boldsymbol{\theta} : d_{\mathrm{H}}(P_0, P_{\boldsymbol{\theta}}) > M\varepsilon_n\}, \ p < p_0 | X_1, \ldots, X_n)$$
$$+ \Pr\left(\left\{\boldsymbol{\theta} : \sum_{j=0}^{N_0} v_j^{2p_0} \theta_j^2 \geq Q, \ d_{\mathrm{H}}(P_0, P_{\boldsymbol{\theta}}) > M\varepsilon_n\right\}, \ p \geq p_0 \,\Big|\, X_1, \ldots, X_n\right)$$
$$+ \Pr\left(\left\{\boldsymbol{\theta} : \sum_{j=0}^{N_0} v_j^{2p_0} \theta_j^2 < Q, \ d_{\mathrm{H}}(P_0, P_{\boldsymbol{\theta}}) > M\varepsilon_n\right\}, \ p \geq p_0 \,\Big|\, X_1, \ldots, X_n\right)$$
$$\leq \Pr(p < p_0 | X_1, \ldots, X_n)$$
$$+ \overline{\pi}_n\left(\left\{\boldsymbol{\theta} : \sum_{j=0}^{N_0} v_j^{2p_0} \theta_j^2 \geq Q\right\} \Big| X_1, \ldots, X_n\right)$$
$$+ \overline{\pi}_n\left(\left\{\boldsymbol{\theta} : \sum_{j=0}^{N_0} v_j^{2p_0} \theta_j^2 < Q, \ d_{\mathrm{H}}(P_0, P_{\boldsymbol{\theta}}) > M\varepsilon_n\right\} \Big| X_1, \ldots, X_n\right)$$
$$\triangleq U_n^{(1)} + U_n^{(2)} + U_n^{(3)}.$$

If $U_n^{(r)} \xrightarrow{\mathrm{P}} 0$ for $r = 1, 2, 3$, then $U_n \xrightarrow{\mathrm{P}} 0$, where all 'in probability' statements are understood to be w.r.t. $P_0^n$. The proof is split into three main steps.

We begin by showing that $U_n^{(1)} \xrightarrow{\mathrm{P}} 0$, namely, that the posterior probability of selecting a model coarser than the best one tends to zero in probability. Note that if $p_0 = 1$, then $U_n^{(1)} = \Pr(p < 1 | X_1, \ldots, X_n) = 0$ a.s. $[P_0^n]$ for all $n \geq 1$. For $p_0 \geq 2$, since $w(m) > 0$ for all $m \in \mathbb{M}$,

$$U_n^{(1)} < \frac{1}{w(0)} \sum_{m<0} w(m) \frac{\int \prod_{i=1}^n f_{\boldsymbol{\theta}_{N_m}}(X_i) \mu_{N_m}(d\boldsymbol{\theta}_{N_m})}{\int \prod_{i=1}^n f_{\boldsymbol{\theta}_{N_0}}(X_i) \mu_{N_0}(d\boldsymbol{\theta}_{N_0})},$$

where the set of integration is understood to be the whole domain. Let
$$R_{m,n} = \frac{\int \prod_{i=1}^{n} f_{\boldsymbol{\theta}_{N_m}}(X_i)\mu_{N_m}(d\boldsymbol{\theta}_{N_m})}{\int \prod_{i=1}^{n} f_{\boldsymbol{\theta}_{N_0}}(X_i)\mu_{N_0}(d\boldsymbol{\theta}_{N_0})}.$$

Since $\mathbb{P}$ is a finite set, for some $m^* < 0$,
$$U_n^{(1)} < \frac{1}{w(0)} \sum_{m<0} w(m) R_{m,n} \leq \frac{\max_{m<0} R_{m,n}}{w(0)} = \frac{R_{m^*,n}}{w(0)} \triangleq \frac{S_n}{w(0)}.$$

It suffices to show that $S_n = o_P(1)$. Using the approximation $e^{n\psi(\boldsymbol{\theta}_{N_m})} \approx e^{n(\theta_0 + \frac{1}{2}\sum_{j=1}^{N_m}\theta_j^2)}$ which is valid for all $m \in \mathbb{M}$ and where $a_n \approx b_n$ means that $a_n/b_n \to 1$ as $n \to \infty$, we obtain that $S_n = T_n + o_P(1)$, with

$$T_n = \prod_{j=1}^{N_0} \left(\frac{n^{-1} + v_j^{-(2p_0+1)}}{n^{-1} + v_j^{-(2p_m+1)}}\right)^{1/2} \exp\left\{\frac{1}{2}\sum_{j=1}^{N_0} b_{j,n}(\bar{\phi}_j)^2\right\}$$
$$\times \prod_{j=N_0+1}^{N_m} (1 + nv_j^{-(2p_m+1)})^{-1/2} \exp\left\{\frac{1}{2}\sum_{j=N_0+1}^{N_m} \frac{(n\bar{\phi}_j)^2}{n + v_j^{2p_m+1}}\right\},$$

where, for simplicity, we have written $m$ instead of $m^*$ and where for $m < 0$,
$$b_{j,n} = n^2[(n + v_j^{2p_m+1})^{-1} - (n + v_j^{2p_0+1})^{-1}] > 0, \qquad 1 \leq j \leq N_0,\ n \geq 1.$$

For later use, note that

(13) $$b_{j,n} < \begin{cases} v_j^{2p_0+1}, & \text{for } n \geq 1, \\ n, & \text{for } 1 \leq j \leq N_0. \end{cases}$$

Recalling the definition of $\bar{\phi}_j$ in Proposition 1, from the inequalities $(\bar{\phi}_j)^2 \leq 2(\bar{\phi}_j - \mathbb{E}_0[\phi_j(X_1)])^2 + 2(\mathbb{E}_0[\phi_j(X_1)])^2$, for all $j \geq 1$, and $x(1+x)^{-1} \leq \log(1 + x) \leq x$, valid for all $x > -1$, it follows that

$$T_n \leq \exp\left\{\frac{1}{2}\sum_{j=1}^{N_0}\left[\frac{v_j^{-(2p_0+1)} - v_j^{-(2p_m+1)}}{n^{-1} + v_j^{-(2p_m+1)}} + 2b_{j,n}(\mathbb{E}_0[\phi_j(X_1)])^2\right]\right\}$$
$$\times \exp\left\{\sum_{j=1}^{N_0} b_{j,n}(\bar{\phi}_j - \mathbb{E}_0[\phi_j(X_1)])^2\right\}$$
$$\times \exp\left\{-\frac{1}{2}\sum_{j=N_0+1}^{N_m}\left[\frac{1}{1 + n^{-1}v_j^{2p_m+1}} - \frac{2n^2(\mathbb{E}_0[\phi_j(X_1)])^2}{n + v_j^{2p_m+1}}\right]\right\}$$
$$\times \exp\left\{\sum_{j=N_0+1}^{N_m} \frac{n^2(\bar{\phi}_j - \mathbb{E}_0[\phi_j(X_1)])^2}{n + v_j^{2p_m+1}}\right\}$$
$$\triangleq T_n^{(1)} \times T_n^{(2)} \times T_n^{(3)} \times T_n^{(4)}.$$

If $\prod_{s=1}^{4} T_n^{(s)} = o_P(1)$, then $T_n = o_P(1)$ and, consequently, $S_n = o_P(1)$. We prove that $T_n^{(1)} = o(1)$. Let $D_0 = \max\{1, \sum_{j=1}^{\infty} v_j^{2p_0} (\mathbb{E}_0[\phi_j(X_1)])^2\}$. Clearly, $D_0 < \infty$ due to (11). Let $n_0$ be the smallest $n$ such that $N_0 \geq 2$. For $n \geq n_0$ and $1 \leq k_n < N_0$ to be specified shortly, recalling (13), we have

$$\sum_{j=1}^{N_0} b_{j,n} (\mathbb{E}_0[\phi_j(X_1)])^2 < v_{k_n} \sum_{j=1}^{k_n} v_j^{2p_0} (\mathbb{E}_0[\phi_j(X_1)])^2$$

$$+ n v_{k_n}^{-2p_0} \sum_{j=k_n+1}^{N_0} v_j^{2p_0} (\mathbb{E}_0[\phi_j(X_1)])^2$$

$$\leq v_{k_n} D_0 + n v_{k_n}^{-2p_0} \sum_{j=k_n+1}^{N_0} v_j^{2p_0} (\mathbb{E}_0[\phi_j(X_1)])^2.$$

Taking $k_n = \lceil (128 D_0)^{-1} n^{1/(2p_0+1)} \rceil$, for $n \geq \max\{n_0, n_1\}$, with $n_1$ the smallest $n$ such that $(128 D_0)^{-1} n^{1/(2p_0+1)} \geq 2$, we have that $v_{k_n} D_0 \leq \frac{1}{64} n^{1/(2p_0+1)}$ and $n v_{k_n}^{-2p_0} \leq (128 D_0)^{2p_0} n^{1/(2p_0+1)}$. For $n \geq \max\{n_0, n_1, n_2\}$, with $n_2$ the smallest $n$ such that $\sum_{j=k_n+1}^{N_0} v_j^{2p_0} (\mathbb{E}_0[\phi_j(X_1)])^2 \leq \frac{1}{64} (128 D_0)^{-2p_0}$, we obtain that $\sum_{j=1}^{N_0} b_{j,n} (\mathbb{E}_0[\phi_j(X_1)])^2 < \frac{1}{32} n^{1/(2p_0+1)}$. Now, note that for $m < 0$,

$$v_j^{-(2p_0+1)} - v_j^{-(2p_m+1)} \leq \begin{cases} 0, & \text{for } j \geq 1, \\ -\frac{1}{2} v_j^{-(2p_m+1)}, & \text{for } j \geq J = \lceil 2^{1/[2(p_0 - p_{-1})]} \rceil \end{cases}$$

so that

$$\frac{1}{2} \sum_{j=1}^{N_0} \frac{v_j^{-(2p_0+1)} - v_j^{-(2p_m+1)}}{n^{-1} + v_j^{-(2p_m+1)}} < -\frac{1}{4} \sum_{j=J}^{\lfloor n^{1/(2p_0+1)} - 1 \rfloor} \frac{v_j^{-(2p_m+1)}}{n^{-1} + v_j^{-(2p_m+1)}}.$$

For $j \leq \lfloor n^{1/(2p_0+1)} - 1 \rfloor$, we have $n^{-1} < v_j^{-(2p_m+1)}$. Also, for $n \geq n_3 = \lceil (2(J+1))^{2p_0+1} \rceil$, we have $J + 1 \leq \frac{1}{2} n^{1/(2p_0+1)}$. Thus, for $n \geq \max\{n_0, n_1, n_2, n_3\}$, combining previous facts, we obtain that

$$0 \leq T_n^{(1)} < \exp\left\{-\frac{1}{4} \sum_{j=J}^{\lfloor n^{1/(2p_0+1)} - 1 \rfloor} \frac{v_j^{-(2p_m+1)}}{n^{-1} + v_j^{-(2p_m+1)}} + \frac{1}{32} n^{1/(2p_0+1)}\right\}$$

$$< \exp\left\{-\frac{1}{8}(\lfloor n^{1/(2p_0+1)} - 1 \rfloor - (J - 1)) + \frac{1}{32} n^{1/(2p_0+1)}\right\}$$

$$\leq \exp\left\{-\frac{1}{32} n^{1/(2p_0+1)}\right\}.$$

Hence, $T_n^{(1)} \to 0$ as $n \to \infty$. We claim that $T_n^{(2)} \xrightarrow{P} 1$. For any $\eta > 0$, by Markov's inequality,

$$P_0^n\left(\sum_{j=1}^{N_0} b_{j,n}(\bar{\phi}_j - \mathbb{E}_0[\phi_j(X_1)])^2 > \eta\right) < \frac{1}{\eta}\sum_{j=1}^{N_0} \frac{b_{j,n}}{n}\mathbb{V}_0[\phi_j(X_1)].$$

By the reverse of Fatou's lemma, the right-hand side goes to zero as $n \to \infty$. To see this, let $\mu$ denote the counting measure on $\mathbb{N}$, endowed with the $\sigma$-field $\mathscr{P}(\mathbb{N})$ of all subsets of $\mathbb{N}$. For each $n \geq 1$, letting

$$f_n(s) = n^{-1} b_{s,n} \mathbb{V}_0[\phi_s(X_1)] I_{\{1,\ldots,N_0\}}(s), \qquad s \geq 1,$$

we can write $\sum_{j=1}^{N_0} n^{-1} b_{j,n} \mathbb{V}_0[\phi_j(X_1)] = \int_{\mathbb{N}} f_n(s)\mu(ds)$. Note that $\{f_n(\cdot), n = 1, 2, \ldots\}$ is a sequence of nonnegative, $\mathscr{P}(\mathbb{N})$-measurable functions such that for every $n \geq 1$,

$$f_n(s) < \mathbb{V}_0[\phi_s(X_1)], \qquad s \geq 1,$$

with $\int_{\mathbb{N}} \mathbb{V}_0[\phi_s(X_1)]\mu(ds) = \sum_{j=1}^{\infty} \mathbb{V}_0[\phi_j(X_1)] < \infty$, due to (12), and

$$\lim_{n \to \infty} f_n(s) = 0, \qquad s \geq 1.$$

Then,

$$\limsup_{n \to \infty} \int_{\mathbb{N}} f_n(s)\mu(ds) \leq \int_{\mathbb{N}} \limsup_{n \to \infty} f_n(s)\mu(ds) = 0.$$

Therefore, $\sum_{j=1}^{N_0} b_{j,n}(\bar{\phi}_j - \mathbb{E}_0[\phi_j(X_1)])^2 \xrightarrow{P} 0$ and by the continuous mapping theorem (CMT), $T_n^{(2)} \xrightarrow{P} 1$. We show that $T_n^{(3)} = o(1)$. Let $n_4$ be the smallest $n$ such that $\sum_{j=N_0+1}^{\infty} v_j^{2p_0}(\mathbb{E}_0[\phi_j(X_1)])^2 \leq 4^{-(p_0+2)}$ and $n_5$ the smallest $n$ such that for $m < 0$, $(N_m/N_0 - 1) \geq 1$. Then for $n \geq \max\{n_3, n_4, n_5\}$,

$$0 \leq T_n^{(3)} < \exp\left\{-\frac{N_m - N_0}{2[1 + 2^{2p_m+1}]} + 4^{-(p_0+2)}n^{1/(2p_0+1)}\right\}$$

$$< \exp\left\{-\frac{1}{4^{p_0+2}}[2N_0 - n^{1/(2p_0+1)}]\right\} \leq \exp\left\{-\frac{n^{1/(2p_0+1)}}{4^{p_0+2}}\right\}.$$

Thus, $T_n^{(3)} \to 0$ as $n \to \infty$. We prove that $T_n^{(4)} = O_P(1)$. For any $\eta > 0$, by Markov's inequality,

$$P_0^n\left(\sum_{j=N_0+1}^{N_m} \frac{n^2(\bar{\phi}_j - \mathbb{E}_0[\phi_j(X_1)])^2}{n + v_j^{2p_m+1}} > \eta\right) < \frac{1}{\eta}\sum_{j=N_0+1}^{N_m} \mathbb{V}_0[\phi_j(X_1)],$$

where the right-hand side goes to zero as $n \to \infty$. By the CMT, $T_n^{(4)} \xrightarrow{P} 1$. Combining all previous results, $T_n \xrightarrow{P} 0$, hence, $U_n^{(1)} \xrightarrow{P} 0$.

The second step consists in showing that for a sufficiently large $Q$, the posterior probability of $\{\boldsymbol{\theta}:\sum_{j=0}^{N_0} v_j^{2p_0}\theta_j^2 \geq Q\}$, under the reduced prior $\overline{\pi}_n$, is asymptotically negligible in probability. Given any $\eta > 0$, by Markov's inequality,

$$P_0^n(U_n^{(2)} > \eta) < \frac{1}{\eta}\mathbb{E}_0^n\left[\overline{\pi}_n\left(\left\{\boldsymbol{\theta}:\sum_{j=0}^{N_0} v_j^{2p_0}\theta_j^2 \geq Q\right\}\Big| X_1,\ldots,X_n\right)\right]$$

$$\leq \frac{1}{\eta Q}\sum_{j=1}^{N_0} v_j^{2p_0}\mathbb{E}_0^n[\mathbb{E}[\theta_j^2|X_1,\ldots,X_n]],$$

where, for $j = 1,\ldots,N_0$,

$$\mathbb{E}[\theta_j^2|X_1,\ldots,X_n] = \sum_{m\geq 0} \overline{w}(m|X_1,\ldots,X_n)\mathbb{E}[\theta_j^2|p=p_m,X_1,\ldots,X_n].$$

Note that conditionally on $p = p_m$, $m \in \mathbb{M}$,

$$\mathbb{E}[\theta_j^2|p=p_m,X_1,\ldots,X_n] \lesssim \frac{1}{n+v_j^{2p_m+1}} + \frac{(n\bar{\phi}_j)^2}{(n+v_j^{2p_m+1})^2}, \quad j=1,\ldots,N_m.$$

Thus, for $j = 1,\ldots,N_0$,

$$\mathbb{E}[\theta_j^2|X_1,\ldots,X_n] \lesssim \frac{1}{n+v_j^{2p_0+1}} + \frac{(n\bar{\phi}_j)^2}{(n+v_j^{2p_0+1})^2}.$$

Since $n^2\mathbb{E}_0^n[(\bar{\phi}_j)^2] = n\mathbb{V}_0[\phi_j(X_1)] + n^2(\mathbb{E}_0[\phi_j(X_1)])^2$, $j \geq 1$, and $\sum_{j=1}^{N_0} v_j^{2p_0}(n+v_j^{2p_0+1})^{-1} < 2^{2p_0+1}$, we have

$$\sum_{j=1}^{N_0} v_j^{2p_0}\left(\frac{1}{n+v_j^{2p_0+1}} + \frac{n^2\mathbb{E}_0^n[(\bar{\phi}_j)^2]}{(n+v_j^{2p_0+1})^2}\right)$$

$$\leq \sum_{j=1}^{N_0}\left(\frac{v_j^{2p_0}(1+\mathbb{V}_0[\phi_j(X_1)])}{n+v_j^{2p_0+1}} + \frac{v_j^{2p_0}n^2(\mathbb{E}_0[\phi_j(X_1)])^2}{(n+v_j^{2p_0+1})^2}\right)$$

$$< \sum_{j=1}^{N_0}\frac{3v_j^{2p_0}}{n+v_j^{2p_0+1}} + \sum_{j=1}^{N_0}v_j^{2p_0}(\mathbb{E}_0[\phi_j(X_1)])^2 < 3\times 2^{2p_0+1} + D_0.$$

Therefore, the probability $P_0^n(U_n^{(2)} > \eta)$ can be made arbitrarily small for all large $n$ by choosing sufficiently large $Q$. Let $Q$ be sufficiently large that $\sum_{j=0}^{\infty} v_j^{2p_0}\theta_{0,j}^2 < Q$. For the same $Q$, define $H_n^c = \{\boldsymbol{\theta}:\sum_{j=0}^{N_0} v_j^{2p_0}\theta_j^2 < Q, d_H(P_0,P_{\boldsymbol{\theta}}) > M\varepsilon_n\}$. In the last step, it remains to be shown that the posterior distribution of $\overline{\pi}_n$ concentrates on $P_0$-centered Hellinger balls at the best rate. Precisely, we prove that $U_n^{(3)} \to 0$, as $n \to \infty$, a.s. $[P_0^\infty]$. Note that

$$U_n^{(3)} < \frac{1}{\overline{w}(0)} \sum_{m \geq 0} \overline{w}(m) \frac{\int_{H_n^c} \prod_{i=1}^n f_{\boldsymbol{\theta}_{N_m}}(X_i) \mu_{N_m}(d\boldsymbol{\theta}_{N_m})}{\int \prod_{i=1}^n f_{\boldsymbol{\theta}_{N_0}}(X_i) \mu_{N_0}(d\boldsymbol{\theta}_{N_0})}.$$

The numerator of the ratio in the summation on the right-hand side of the above inequality can be bounded above using condition (16), as in the proof of Theorem 1. To bound the denominator below, we can use the same arguments as in the proof of Theorem 3, replacing $N$ with $N_0 = \lceil n^{1/(2p_0+1)} \rceil$ and taking $\bar{B}_1^2 < \min\{B_1^2/2, (\sqrt{Q} - \sqrt{Q - \delta_0})^2/16^{p_0}\}$. Then for $n$ sufficiently large that $\sum_{j=N_0+1}^\infty v_j^{2p_0} \theta_{0,j}^2 \leq B_1^2/2$, we have

$$\mu_{N_0}\left(\left\{\boldsymbol{\theta}_{N_0} \in E_{p_0, N_0}(Q) : \sum_{j=1}^\infty (\theta_j - \theta_{0,j})^2 \leq B_1^2 \varepsilon_n^2\right\}\right)$$
$$\geq \mu_{N_0}\left(\left\{\boldsymbol{\theta}_{N_0} : \sum_{j=1}^{N_0} (\theta_j - \theta_{0,j})^2 \leq \bar{B}_1^2 \varepsilon_n^2\right\}\right) \gtrsim e^{-c_1 n \varepsilon_n^2},$$

with $c_1$ depending on $\boldsymbol{\theta}_0$. Thus, for a suitable constant $c > 0$, $U_n^{(3)} \lesssim e^{-cn\varepsilon_n^2}$ for all but finitely many $n$ along almost all sample paths when sampling from $P_0$. This completes the proof. □

REMARK 6. If $f_0$ simultaneously has the following series expansions

$$f_0(x) = \boldsymbol{\beta}_0 \cdot \boldsymbol{\phi}(x) = \exp\{\boldsymbol{\theta}_0 \cdot \boldsymbol{\phi}(x) - \psi(\boldsymbol{\theta}_0)\}, \qquad x \in [0,1],$$

where $\boldsymbol{\beta}_0 = (\beta_{0,0}, \beta_{0,1}, \ldots)$ has coordinates $\beta_{0,0} = 1$ and $\beta_{0,j} = \mathbb{E}_0[\phi_j(X_1)]$ for $j \geq 1$, then condition (11) implies that $\boldsymbol{\beta}_0$ lies in $E_{p_0}$.

REMARK 7. Since a finite set $\mathbb{P}$ of possible values for $p$ is considered, the choice of weights $w(m)$ is not relevant. In the present asymptotic setting, any set of positive weights achieves the same result.

Since the posterior distribution does not converge exponentially fast in Theorem 5, the rate of convergence for the posterior expected density cannot be derived as easily as in the previous cases. We therefore resort to another estimator that is Bayesian in the sense that it is based on the posterior distribution. The following construction closely follows that in [4], pages 544–545. For a positive sequence $\delta_n \to 0$ as $n \to \infty$, let $H_{\delta_n}(P_{\boldsymbol{\theta}}) = \{P_{\boldsymbol{\theta}'} : d_{\mathrm{H}}(P_{\boldsymbol{\theta}}, P_{\boldsymbol{\theta}'}) \leq \delta_n\}$ and define

$$\delta_n^* = \inf\{\delta_n : \Pi_n(H_{\delta_n}(P_{\boldsymbol{\theta}})|X_1, \ldots, X_n) \geq 3/4 \text{ for some } P_{\boldsymbol{\theta}}\}.$$

Take any $\hat{P}_n$ satisfying the following condition:

$$\Pi_n(H_{\delta_n^* + n^{-1}}(\hat{P}_n)|X_1, \ldots, X_n) \geq 3/4.$$

As subsequently stated, such an estimator, whose definition does not require knowledge of $p_0$, attains the optimal pointwise rate of convergence $n^{-p_0/(2p_0+1)}$, adapting to the unknown smoothness of the true density.

COROLLARY 3. *If the conditions of Theorem 5 are satisfied, then for a sufficiently large constant $M > 0$, $P_0^n(d_{\mathrm{H}}(P_0, \hat{P}_n) > Mn^{-p_0/(2p_0+1)}) \to 0$ as $n \to \infty$.*

PROOF. See the proof of Theorem 4 in [4], page 545. □

**6. Closing remarks.** This paper focuses on the estimation of densities in periodic Sobolev classes. The problem is approached through the use of an orthonormal series expansion for the log-density with single priors on the coefficients. The posterior expected density is shown to attain the optimal minimax rate of convergence under Hellinger loss for several priors.

As mentioned in Remark 1, an interesting finding of the paper is that a sample-size-dependent direct Gaussian prior leads to a Bayes' estimator achieving the optimal minimax rate in this problem, in contrast to the infinitely many normal means problem investigated by Zhao [26], who has shown that there is no Gaussian prior supported on $E_p$ such that the corresponding Bayes' estimator attains the optimal minimax rate. Optimality for the Bayes' density estimator follows from uniform exponential convergence of the posterior distribution over suitable ellipsoids. In the infinitely many normal means problem, the rate of convergence for the Bayes' estimator is derived directly from the study of the risk function and uniformity holds over any $E_p(Q)$ provided the power of the prior variances exactly matches the assumed degree of smoothness, which is not the case if the prior is supported on $E_p$.

Another interesting result concerns adaptation. We have shown that the posterior distribution of a sample-size-dependent prior achieves the best pointwise rate $n^{-p_0/(2p_0+1)}$, regardless of the value of $p_0 \in \mathbb{P}$, for every $\boldsymbol{\theta}_0 \in E_{p_0}$ satisfying conditions (11) and (12). In a recent paper, Huang [16] has obtained results on posterior rates of convergence for density estimation using the method of exponentials, with priors on the coefficients of the log-density expansion via wavelets, the coefficients lying in a Besov space $B_{2,2}^\alpha$ with $\alpha \in (0, 1)$. This method is suitable for estimating spatially inhomogeneous density functions, while we consider smooth, periodic functions. Huang does not put a prior on $\alpha$, instead she constructs a sieve prior with mixing parameter given by the dimension of the exponential family and the ball radius. Even though the rate she obtains has an extra $(\log n)^{1/2}$-factor, her result is valid for all points in $B_{2,2}^\alpha$. Our result, although achieving a better rate, is restricted to points in $E_{p_0}$ also satisfying the aforementioned conditions.

## APPENDIX

LEMMA A.1. *For any pair $\boldsymbol{\theta}$, $\boldsymbol{\theta}' \in E_p(Q)$,*

$$(14) \qquad K(P_{\boldsymbol{\theta}'} \| P_{\boldsymbol{\theta}}) < \frac{1}{2} e^{4B} \sum_{j=1}^{\infty} (\theta_j' - \theta_j)^2.$$

*Consequently,*

$$(15) \qquad d_{\mathrm{H}}^2(P_{\boldsymbol{\theta}'}, P_{\boldsymbol{\theta}}) \| f_{\boldsymbol{\theta}'} / f_{\boldsymbol{\theta}} \|_{\infty} < \frac{1}{2} e^{8B} \sum_{j=1}^{\infty} (\theta_j' - \theta_j)^2.$$

PROOF. We use inequality (3.2) from Lemma 1 of Barron and Sheu [3], pages 1355–1356. If $\| \log(f_{\boldsymbol{\theta}'} / f_{\boldsymbol{\theta}}) \|_{\infty} < \infty$, then for any constant $c$,

$$K(P_{\boldsymbol{\theta}'} \| P_{\boldsymbol{\theta}}) \leq \frac{1}{2} e^{\| \log(f_{\boldsymbol{\theta}'}/f_{\boldsymbol{\theta}}) - c \|_{\infty}} \int_0^1 f_{\boldsymbol{\theta}'}(x) \left( \log \frac{f_{\boldsymbol{\theta}'}(x)}{f_{\boldsymbol{\theta}}(x)} - c \right)^2 dx.$$

Note that for any pair $\boldsymbol{\theta}$, $\boldsymbol{\theta}' \in E_p(Q)$,

$$|[\psi(\boldsymbol{\theta}) - \theta_0] - [\psi(\boldsymbol{\theta}') - \theta_0']| \leq \|(\boldsymbol{\theta} - \boldsymbol{\theta}') \cdot \boldsymbol{\phi} - (\theta_0 - \theta_0')\|_{\infty} < 2B,$$

thus,

$$\left\| \log \frac{f_{\boldsymbol{\theta}'}}{f_{\boldsymbol{\theta}}} \right\|_{\infty} = \|[(\boldsymbol{\theta}' - \boldsymbol{\theta}) \cdot \boldsymbol{\phi}] + [\psi(\boldsymbol{\theta}) - \psi(\boldsymbol{\theta}')]\|_{\infty} < 4B < \infty.$$

Take $c = [\psi(\boldsymbol{\theta}) - \theta_0] - [\psi(\boldsymbol{\theta}') - \theta_0']$. Using the fact that $\sup_{\boldsymbol{\theta} \in E_p(Q)} \|f_{\boldsymbol{\theta}}\|_{\infty} < e^{2B}$ and Parseval's relation, we obtain (14). Obviously, the same bound holds true for $K(P_{\boldsymbol{\theta}} \| P_{\boldsymbol{\theta}'})$. A similar remark applies to inequality (15), which follows from (14) because $d_{\mathrm{H}}^2(P_{\boldsymbol{\theta}'}, P_{\boldsymbol{\theta}}) \leq K(P_{\boldsymbol{\theta}'} \| P_{\boldsymbol{\theta}})$ and $\|f_{\boldsymbol{\theta}'}/f_{\boldsymbol{\theta}}\|_{\infty} < e^{4B}$. □

Theorem A.1 below is an almost sure version of Theorem 2.1 in [15], page 1239 (see also Theorem 2.2 in [16], page 505). Before stating the theorem, we recall that if $(S, d)$ is a semi-metric space and $C$ a totally bounded subset of $S$, then for any given $\varepsilon > 0$, the $\varepsilon$-*packing number* of $C$, denoted by $D(\varepsilon, C, d)$, is defined as the largest integer $m$ such that there exists a set $\{s_1, \ldots, s_m\} \subseteq C$ with $d(s_k, s_l) > \varepsilon$ for all $k, l = 1, \ldots, m$, $k \neq l$. The $\varepsilon$-*capacity* of $(C, d)$ is defined as $\log D(\varepsilon, C, d)$.

THEOREM A.1. *Let $\mathscr{P}$ be a class of probability measures that possess densities relative to some $\sigma$-finite reference measure $\nu$ on a sample space $(\mathscr{X}, \mathscr{A})$. Let $d$ stand for either the $L_1$- or the Hellinger metric on $\mathscr{P}$. Let $\Pi_n$ be a sequence of priors on $(\mathscr{P}, \mathscr{B})$, where $\mathscr{B}$ is the Borel $\sigma$-field on $\mathscr{P}$. For $P_0 \in \mathscr{P}$, let $f_0$ denote its density. Suppose that for positive sequences* $\bar{\varepsilon}_n, \tilde{\varepsilon}_n \to 0$ with $n \min\{\bar{\varepsilon}_n^2, \tilde{\varepsilon}_n^2\} \to \infty$ and $\sum_{n=1}^{\infty} \exp(-En\tilde{\varepsilon}_n^2) < \infty$ for every $E > 0$, constants $c_1, c_2, c_3, c_4 > 0$ and sets $\mathscr{P}_n \subseteq \mathscr{P}$, we have

$$\log D(\bar{\varepsilon}_n, \mathscr{P}_n, d) \leq c_1 n \bar{\varepsilon}_n^2, \tag{16}$$

$$\Pi_n(\mathscr{P} \setminus \mathscr{P}_n) \leq c_2 e^{-(c_3+4)n\tilde{\varepsilon}_n^2}, \tag{17}$$

$$\Pi_n(N(P_0; \tilde{\varepsilon}_n^2)) \geq c_4 e^{-c_3 n \tilde{\varepsilon}_n^2}, \tag{18}$$

where $N(P_0; \tilde{\varepsilon}_n^2) = \{P : d_{\mathrm{H}}^2(P_0, P) \|f_0/f_P\|_\infty \leq \tilde{\varepsilon}_n^2\}$ with $f_P = dP/d\nu$. Then, for $\varepsilon_n = \max\{\bar{\varepsilon}_n, \tilde{\varepsilon}_n\}$ and a sufficiently large constant $M > 0$, the posterior probability

$$\Pi_n(\{P : d(P_0, P) > M\varepsilon_n\} | X_1, \ldots, X_n) \to 0 \qquad as\ n \to \infty,$$

$P_0^\infty$-almost surely.

PROOF OF THEOREM 1. We appeal to Theorem A.1 and show that the conditions listed earlier are satisfied for $\bar{\varepsilon}_n = \tilde{\varepsilon}_n = \varepsilon_n = n^{-p/(2p+1)}$. Condition (16) is verified for $\mathscr{P}_n = \mathscr{F}$. It is easily seen that for some constant $K > 0$ depending only on $p$ and $L$, $\int_0^1 (f_{\boldsymbol{\theta}}^{(p)}(x))^2\, dx < K^2$ for all $\boldsymbol{\theta} \in E_p(Q)$. Besides, for any pair $P_{\boldsymbol{\theta}'}, P_{\boldsymbol{\theta}} \in \mathscr{F}$ such that $d_{\mathrm{H}}(P_{\boldsymbol{\theta}'}, P_{\boldsymbol{\theta}}) > \varepsilon_n$, a simple calculation shows that $\|f_{\boldsymbol{\theta}'} - f_{\boldsymbol{\theta}}\|_\infty \geq \|f_{\boldsymbol{\theta}'} - f_{\boldsymbol{\theta}}\|_2 > 2e^{-B}\varepsilon_n$; see [5], page 252, for the monotone convergence of the $L_q$-norm, $q \geq 1$, to the essential supremum norm w.r.t. $\lambda$ on $[0,1]$, $\|\cdot\|_q \uparrow \|\cdot\|_{L_\infty}$. Then by a result due to Birman and Solomjak [6] (see also [20], pages 22–23), for a suitable constant $c > 0$,

$$\log D(\varepsilon_n, \mathscr{F}, d_{\mathrm{H}}) \leq \log D(2e^{-B}\varepsilon_n, \mathscr{F}, \|\cdot\|_\infty) \leq c\varepsilon_n^{-1/p} = cn\varepsilon_n^2.$$

Condition (17) is trivially verified. Finally, recalling that $B_1^2 = e^{-8B}$, condition (18) follows from

$$c_2 e^{-c_1 n \varepsilon_n^2} \leq \pi_n\left(\left\{\boldsymbol{\theta} : \sum_{j=1}^\infty (\theta_j - \theta_{0,j})^2 \leq e^{-8B}\varepsilon_n^2\right\}\right)$$

$$\leq \pi_n(\{\boldsymbol{\theta} : d_{\mathrm{H}}^2(P_0, P_{\boldsymbol{\theta}}) \|f_0/f_{\boldsymbol{\theta}}\|_\infty \leq \varepsilon_n^2\}) = \Pi_n(N(P_0; \varepsilon_n^2)),$$

where (1) and (15) have been applied. □

**Acknowledgments.** The author would like to thank two anonymous referees, an Associate Editor and the Editor for providing constructive comments and insightful suggestions.

ISTITUTO DI METODI QUANTITATIVI
UNIVERSITÀ "L. BOCCONI," MILANO
VIALE ISONZO 25
I-20135 MILANO
ITALY
E-MAIL: catia.scricciolo@uni-bocconi.it